# PARAMETRIC MODELING APPROACH TO COVID-19 PANDEMIC DATA

N. I. Badmus, O. Faweya and S. A. Ige

Department of Mathematics, University of Lagos, Akoka, Nigeria
Department of Statistics, Ekiti State University, Ado-Ekiti, Ekiti State, Nigeria
Department of Mathematics, Yaba College of Technology, Lagos, Nigeria
Correspondence Email: nibadmus@unilag.edu.ng

Abstract

The problem of skewness is common among clinical trials and survival data which has being the research focus derivation and proposition of different flexible distributions. Thus, a new distribution called Extended Rayleigh Lomax distribution is constructed from Rayleigh Lomax distribution to capture the excessiveness of some survival data. We derive the new distribution by using beta logit function proposed by Jones (2004). Some statistical properties of the distribution such as probability density function, cumulative density function, reliability rate, hazard rate, reverse hazard rate, moment generating functions, likelihood functions, skewness, kurtosis and coefficient of variation are obtained. We also performed the expected estimation of model parameters by maximum likelihood; goodness of fit and model selection criteria including Anderson Darling (AD), CramerVon Misses (CVM), Kolmogorov Smirnov (KS), Akaike Information Criterion (AIC), Bayesian Information Criterion (BIC) and Consistent Akaike Information Criterion (CAIC) are employed to select the better distribution from those models considered in the work. The results from the statistics criteria show that the proposed distribution performs better with better representation of the States in Nigeria COVID-19 death cases data than other competing models.

## 1. Introduction

In survival analysis, problems are encountered in the analysis of clinical data because distributions proposed are not flexible enough to follow the movement of the data to give accurate results. In the light of this, there is need to develop a more flexible parametric model using COVID-19 data for example. In recent times, there was outbreak of the third wave of COVID-19 pandemic called Delta Variant after the second wave generating a global outcry. Many researches/ works have been done by several researchers since the breakup of the pandemic in December 2019 from various fields such as: Medicine, Statistics, Economics etc., with different ideas, models, methods and approaches in their respective works. These include: Badmus *et al.* 2020, Dey *et al.* 2020, WHO, 2020, Yoo, 2020 amongst others. Most clinical data are always skewed, thus a new distribution is constructed and generated from a parent distribution called Rayleigh Lomax (RL) distribution by Kawsar *et al*. (2018) is generated using beta link function introduced by Jones (2004). This is expected to have different shapes for the survival and hazard rate functions. More parameters are added to the parent distribution, the flexibility and the ability of the distribution to model real life data are established

## 2. Material and Methods

There are several methods in literature which have been used by many researchers. In this study, we consider beta logit function introduced by Jones (2004), which can jointly convolute two or more distributions.

### 2.1 Properties of Extended Rayleigh Lomax (ERL) Distribution

### 2.1 Probability Density Function

The probability density function of the above distribution is obtained using the beta link function given as:

$$g(x, a, b) = \frac{(k(x)[K(x)]^{(a-1)}[1-K(x)]^{(b-1)})}{B(a,b)} \quad (1)$$

$$x, \alpha, b, \theta, \lambda, \beta \simeq ERL(x, \alpha, b, \theta, \lambda, \beta > 0)$$

$$K(x) = 1 - e^{-\frac{\beta}{2}\left(\frac{\theta}{\theta+x}\right)^{-2\lambda}} \text{ and } k(x) = \frac{\beta\lambda}{\theta}\left(\frac{\theta}{\theta+x}\right)^{-2\lambda+1} e^{-\frac{\beta}{2}\left(\frac{\theta}{\theta+x}\right)^{-2\lambda}} \quad x \geq -\theta \text{ and } \theta, \lambda, \beta > 0$$

where $K(x)$ and $k(x)$ are the cdf and pdf of the parent distribution respectively, $a$ and $b$ are additional shape parameters to the parent distribution.

$$g(x, \alpha, b, \theta, \lambda, \beta) = \frac{1}{B(a,b)}[K(x)]^{a-1}(1-K(x))^{b-1}k(x) < x, \alpha, b, \theta, \lambda, \beta > 0 \quad (2)$$

where θ, λ and β are scale and shape parameters while a and b are the new shape parameters introduced to the distribution. Then (2) becomes extended Rayleigh Lomax (ERL) Distribution

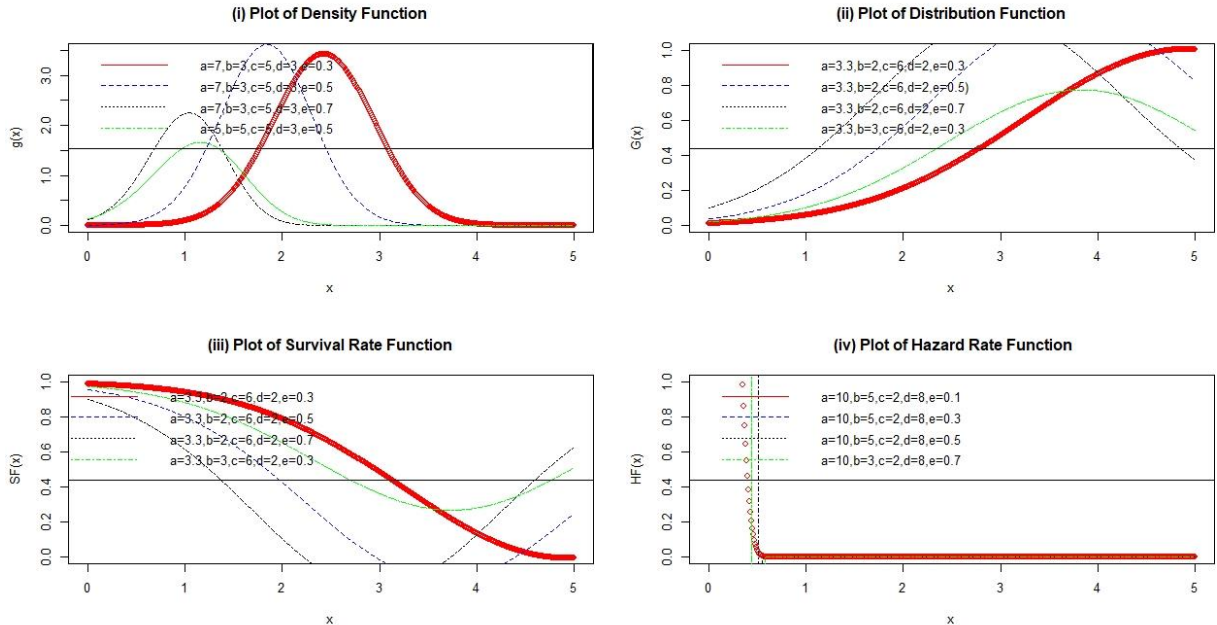

Figure 1 depicts: (i) the pdf ,(ii) the cdf , (iii) the survival and (iv) hazard function in plots of the ERL distribution.

Some new distributions from the extended Rayleigh Lomax distribution are mentioned below:

(a) When $a = 1$ in (2), we have the Lehmann Type II Rayleigh Lomax Distribution.

(b) If $b = 1$ in (2), we get the pdf of exponentiated Rayleigh Lomax distribution

(c) When β=1 in (2), this consists of beta Lomax distribution

(d) When $\frac{\beta}{\theta} = \beta, \lambda = 1 \; and \; \frac{\theta+x}{\theta} = 2x$ in (2), it then becomes beta Rayleigh distribution

(e) If $a = b = 1$ in (2), it yields Rayleigh Lomax distribution which is the parent distribution. (see Kawsar *et al.* (2018)).

(f) When $b = \beta = 1$ in (2), this consists exponential Lomax distribution by El-Bassiouny *et al.* (2015)

(g) If $a = b = \lambda = \theta = 1$ in (2), it yields Rayleigh Lomax distribution which is the parent distribution. (see Siddiqui, 1962).

## 2.2 The Cumulative Distribution (CDF) of BRL Distribution

The associative cumulative distribution function cdf in (2) is given as

$$P(X \leq x) = \int_0^x f(k)dk = G(x, a, b, \theta, \lambda, \beta)$$

$$= \int_0^x \frac{1}{B(a,b)} [K(x)]^{a-1} (1 - K(x))^{b-1} k(x) dk \qquad (3)$$

We set

$$k(x) = 1 - e^{-\frac{\beta}{2}\left(\frac{\theta}{\theta+x}\right)^{-2\lambda}} \qquad (4)$$

$$\frac{dk}{dx} = \frac{\beta\lambda}{\theta} \left(\frac{\theta}{\theta+x}\right)^{-2\lambda+1} e^{-\frac{\beta}{2}\left(\frac{\theta}{\theta+x}\right)^{-2\lambda}} \left[1 - e^{-\frac{\beta}{2}\left(\frac{\theta}{\theta+x}\right)^{-2\lambda}}\right]$$

$$\frac{dk}{dx} = \frac{\beta\lambda}{\theta} \left(\frac{\theta}{\theta+x}\right)^{-2\lambda+1} e^{-\frac{\beta}{2}\left(\frac{\theta}{\theta+x}\right)^{-2\lambda}} - \left[e^{-\frac{\beta}{2}\left(\frac{\theta}{\theta+x}\right)^{-2\lambda}}\right]^2$$

Also

$$,dx = \frac{dk}{\frac{\beta\lambda}{\theta}\left(\frac{\theta}{\theta+x}\right)^{-2\lambda+1} e^{-\frac{\beta}{2}\left(\frac{\theta}{\theta+x}\right)^{-2\lambda}} \left[1 - e^{-\frac{\beta}{2}\left(\frac{\theta}{\theta+x}\right)^{-2\lambda}}\right]}$$

Putting dx in equation (2), we realize:

$$g(x, \alpha, b, \theta, \lambda, \beta) = \frac{1}{B(a,b)} \int_0^\infty \left[1 - e^{-\frac{\beta}{2}\left(\frac{\theta}{\theta+x}\right)^{-2\lambda}}\right]^{a-1} \left[e^{-\frac{\beta}{2}\left(\frac{\theta}{\theta+x}\right)^{-2\lambda}}\right]^{b-1} dk \qquad (5)$$

And k in equation (4) becomes

$$k(x) = \frac{dk(x)}{dx} = \frac{\beta\lambda}{\theta} \left(\frac{\theta}{\theta+x}\right)^{-2\lambda+1} e^{-\frac{\beta}{2}\left(\frac{\theta}{\theta+x}\right)^{-2\lambda}}$$

Equation (5) can be expressed as

$$g(x, \alpha, b, \theta, \lambda, \beta) = \frac{1}{B(a,b)} K^{a-1} (1-K)^{b-1} \frac{dk}{dx} \qquad (6)$$

Now, putting (2) in (6), we get

$$G(x, a, b, \theta, \lambda, \beta) = P(X \leq x) = \int_0^x \frac{1}{B(a,b)} K^{a-1} (1-K)^{b-1} \frac{dk}{dx}$$

$$= \int_0^x \frac{1}{B(a,b)} K^{a-1}(1-K)^{b-1} \frac{dk}{dx}$$

Where $B(x,a,b) = \int_0^x K^{a-1}(1-K)^{b-1} dk$ and it is called the incomplete beta function.

$$G(x,a,b,\theta,\lambda,\beta) = \int_0^x \frac{K^{a-1}(1-K)^{b-1} dk}{B(a,b)} = \frac{B(k,a,b)}{B(a,b)} \tag{7}$$

Expression (7) becomes the cumulative distribution function of ERL distribution.

## 2.3 The Survival Rate/Reliability Function

The reliability function of BRL distribution is given by

$$Sur(x,a,b,\theta,\lambda,\beta) = 1 - G(x,a,b,\theta,\lambda,\beta) = 1 - \int_0^x f(k) dk$$

$$1 - \int_0^x \frac{K^{a-1}(1-K)^{b-1} dk}{B(a,b)} = 1 - \frac{B(k,a,b)}{B(a,b)}$$

$$Sur(x,a,b,\theta,\lambda,\beta) = \frac{B(a,b) - B(x,a,b)}{B(a,b)} \tag{8}$$

## 2.4 The Hazard Rate Function

$$hez(x,a,b,\theta,\lambda,\beta) = \frac{g(x,a,b,\theta,\lambda,\beta)}{1 - G(x,a,b,\theta,\lambda,\beta)} = \frac{\frac{K^{a-1}(1-K)^{b-1} k}{B(a,b)}}{\frac{B(a,b) - B(x,a,b)}{B(a,b)}}$$

$$= \frac{K^{a-1}(1-K)^{b-1} k}{B(a,b) - B(x,a,b)} \tag{9}$$

## 2.5 The Reversed Hazard Rate Function

$$Rhaz(x,a,b,\theta,\lambda,\beta) = \frac{g(x,a,b,\theta,\lambda,\beta)}{G(x,a,b,\theta,\lambda,\beta)} = \frac{\frac{K^{a-1}(1-K)^{b-1} k}{B(k,a,b)}}{B(a,b)}.$$

$$= \frac{B(a,b) K^{a-1}(1-K)^{b-1} k}{B(k,a,b)} \tag{10}$$

## 2.6 Testing the Trueness of the PDF of ERL Distribution

The ERL distribution is a probability density function with the use of:

$$\int_0^\infty g_{ERL}(y)dy = 1 \qquad (11)$$

Jones (2004) in his generalized beta distribution of first kind is given by:

$$g_x(x, a, b, u) = [B(a,b)]^{-u}[K(y)]^{au-1}[1-K(y)^u]^{s-1}K(y) \quad 0 < y < 1$$

where $a, b$ and $u > 0$, therefore differentiating $(a)$ above, we obtain

$$g_{ERLD}(y) = [B(a,b)]^{-u} G^{au-1}(1-G^u)^{b-1}\frac{dG}{dy}$$

$$\int_{-\infty}^{\infty} g_{ERLD}(y) = \int_{-\infty}^{\infty} \frac{u}{B(a,b)} G^{au-1}(1-G^u)^{b-1} dG$$

Putting $M = G^u$, then differentiating M with respect to G

$$\frac{dM}{dG} = UG^{u-1}$$

$$dG = \frac{dM}{UG^{u-1}}$$

$$G = M^{\frac{1}{u}}$$

$$\int_{-\infty}^{\infty} g_{ERLD}(y)dy = \int_0^1 \frac{u}{B(a,b)} \left(M^{\frac{1}{u}}\right)^{au-1}(1-M)^{b-1}\frac{dM}{UG^{u-1}}$$

$$= [B(a,b)]^{-1} \int_0^1 \frac{M^{a-\frac{1}{u}}(1-M)^{b-1}}{M^{1-\frac{1}{u}}} dM$$

$$= [B(a,b)]^{-1} \int_0^1 M^{a-1}(1-M)^{b-1} dM$$

$$\int_0^1 M^{a-1}(1-M)^{b-1} dM = B(a,b)$$

Therefore, $$g_{ERLD}(y) = \frac{B(a,b)}{B(a,b)} = 1$$

Hence, the $g_{ERL}$ distribution has a true continuous probability density function.

## 3. Moments and Generating Function

In this section, we derive and obtain the moment generating function (mgf) of the distribution $m(t) = E(e^{ty})$ and the general $rth$ moment of a beta generated distribution defined by Hosking (1990)

$$\mu_r^1 = B(a,b)^{-1} \int_0^1 [F^{-1}(y)]^r y^{a-1} [1-y]^{b-1} dy \qquad (12)$$

Cordeiro et al. (2011) also discussed another mgf for generated beta distribution.

$$m(t) = B(a,b)^{-1} \sum_{j=0}^{\infty} (-1)^j \binom{b-1}{j} e(q,r; aj-1) \qquad (13)$$

where,

$$e(q,r) = \int_{-\infty}^{\infty} e^{ty} [F(y)]^m f(y) dy$$

then,

$$M_y(t) = B(a,b)^{-1} \sum_{j=0}^{\infty} (-1)^j \binom{b-1}{j} \int_{-\infty}^{\infty} e^{ty} [F(y)]^{a(j+1)-1} f(y) dy \qquad (14)$$

Putting the pdf and cdf of the Extended Rayleigh Lomax distribution into equation (14), we get

$$M_{ERLD(y)}(t) = B(a,b)^{-1} \sum_{j=0}^{n} (-1)^j \binom{b-1}{j} \int_{-\infty}^{\infty} e^{ty} [F(y)]^{a(j+1)-1} f(y) \qquad (15)$$

If $a = b = 1$ in equation (14) that becomes the moment generating function of the baseline distribution.

Hence, the $rth$ moment of the ERL distribution is obtained, since the moment generating function of the parent distribution is given by

$$M_y(t) = \sum_{j=0}^{\infty} \frac{t^j}{j!} \int_0^{\infty} y^j f(y, \beta, \lambda, \theta) dy \sum_{c=0}^{\infty} \frac{t^c}{c!} \sum_{j=0}^{\infty} (bC_j) \theta^j \left(\frac{2}{\beta}\right)^{\frac{1}{2\lambda}} (-\theta)^{b-j} r\left(\frac{j}{2\lambda} + 1\right) \qquad (16)$$

Equation (16) can be re written as

$$M_{ERLD(y)}(t) = B(a,b)^{-1} \sum_{j=0}^{n} (-1)^j \binom{b-1}{j} \int_{-\infty}^{\infty} e^{ty} [F(y)]^{a(j+1)-1} \sum_{i=0}^{\infty} \frac{t^i}{i!} \sum_{h=0}^{\infty} (iC_h) \theta^h \left(\frac{2}{\beta}\right)^{\frac{h}{2\lambda}} (-\theta)^{i-h} r\left(\frac{h}{2\lambda} + 1\right)$$

$$= B(a,b)^{-1} \sum_{i=0}^{\infty} \sum_{j=0}^{\infty} \sum_{h=0}^{\infty} (-1)^j \binom{b-1}{i} \frac{t^i}{i!} (iC_h) \theta^h \left(\frac{2}{\beta}\right)^{\frac{h}{2\lambda}} (-\theta)^{i-h} r\left(\frac{h}{2\lambda} + 1\right) \cdot [F(y)]^{a(j+1)-1} \quad (17)$$

and the $rth$ moment of ERL distribution is obtained from equation (17)

$$\mu^i_{ERLD(r)} = E(y^r) = B(a,b)^{-1} \sum_{j=0}^{\infty} \sum_{h=0}^{\infty} (-1)^j \binom{b-1}{i} \cdot [F(y)]^{a(j+1)-1} \frac{t^i}{r!} \binom{r}{h} \theta^h \left(\frac{2}{\beta}\right)^{\frac{h}{2\lambda}} (-\theta)^{r-h} r\left(\frac{h}{2\lambda} + 1\right) \quad (18)$$

Letting $a = b = 1$ in (18) gives the rth moment of the baseline distribution by Kawsar *et al.* (2018)

$$\mu^i_r = E(y^r) = \sum_{h=0}^{\infty} \binom{r}{h} \theta^h \left(\frac{2}{\beta}\right)^{\frac{h}{2\lambda}} (-\theta)^{r-h} r\left(\frac{h}{2\lambda} + 1\right)$$

Other measures such as the Skewness (SK$_{ERLD}$)$(a, b, \beta, \lambda, \theta)$ and Kurtosis (KT$_{ERLD}$) $(a, b, \beta, \lambda, \theta)$ are also obtained below:

The $rth$ moment of the ERL distribution is written as:

$$\mu^i_{ERLD(r)} = \int_0^{\infty} y^r F_{ERLD}(y) dy$$

That is,

$$\mu^i_{ERLD(r)} = \int_0^{\infty} y^r \left([B(a,b)]^{-1} (K(Y))^{a-1} (1 - K(y))^{b-1} dk(y)\right)$$

where, $\quad K(y) = (1 - w(y))^{-u}$

i.e $\quad W(y) = e^{-\frac{\beta}{2}\left(\frac{\theta}{\theta+y}\right)}$ and $u = 2\lambda$

therefore,

$$\mu^i_{ERLD(r)} = \frac{(rC_h)\theta^h \left(\frac{2}{\beta}\right)^{\frac{h}{u}} (-\theta)^{r-h} r(\frac{h}{u}+1)}{B(a,b)} \sum_{j=0}^{\infty} \sum_{h=0}^{\infty} (-1)^j \binom{b-1}{i} \left([(1-w(y))^{-u}]^{a(j+1)-1}\right)$$

$$= Z\left(\binom{r}{h} \theta^h \left(\frac{2}{\beta}\right)^{\frac{h}{u}} (-\theta)^{r-h} r(\frac{h}{u}+1)\right) \quad (19)$$

Where,

$$Z = \frac{\sum_{j=0}^{\infty} \sum_{h=0}^{\infty} (-1)^j \binom{b-1}{i} \left([(1-w(y))^{-u}]^{a(j+1)-1}\right)}{B(a,b)}$$

At the same time, the first four central moments $\mu^i_r = 1, 2, 3, 4$ are obtained through (17) as:

Furthermore, the mean and second to fourth moments of the ERL distribution are given as follows:

$\mu = \mu_1^|$, $\mu_2 = \mu_2^| - \mu^2$, $\mu_3 = \mu_3^| - 3\mu\mu_2^| + 2\mu^3$, and $\mu_4 = \mu_4^| - 4\mu\mu_3^| + 6\mu^2\mu_2^| - 3\mu^4$

$$\mu_1^| = Z\left(\binom{1}{h}\theta^h \left(\frac{2}{\beta}\right)^{\frac{h}{u}} (-\theta)^{1-h} r(\frac{h}{u}+1)\right) \quad (20)$$

$$\mu_2^| = Z\left(\binom{2}{h}\theta^h \left(\frac{2}{\beta}\right)^{\frac{h}{u}} (-\theta)^{2-h} r(\frac{h}{u}+1)\right) \quad (21)$$

$$\mu_3^| = Z\left(\binom{3}{h}\theta^h \left(\frac{2}{\beta}\right)^{\frac{h}{u}} (-\theta)^{3-h} r\left(\frac{h}{u}+1\right)\right) \quad (22)$$

$$\mu_4^| = Z\left(\binom{4}{h}\theta^h \left(\frac{2}{\beta}\right)^{\frac{h}{u}} (-\theta)^{4-h} r(\frac{h}{u}+1)\right) \quad (23)$$

Other measures such as skewness, kurtosis and coefficient of variation of the ERL distribution are given below:

### 3.8 Skewness of the ERL Distribution

The skewness is a means of measuring non symmetry of the distribution. The skewness is given by:

$$SK_{ERLD} = \frac{\mu_3}{\mu_2^{1.5}} \quad (24)$$

### 3.9 Kurtosis of the ERL Distribution

The kurtosis is another measure that measures the peak of the distribution. The kurtosis of the BRL distribution is given as:

$$KT_{ERLD} = \frac{\mu_4}{\mu_2^2} - 3 \quad (25)$$

### 3.10 Coefficient of Variation of the ERL Distribution

This is also a measure of variability of a probability distribution. The CV of the ERL distribution is given as:

$$CV_{ERLD} = \frac{\sqrt{\mu_2}}{\mu} \quad (26)$$

## 4. Estimation of Parameter

We made attempt to derive the maximum likelihood estimates (MLEs) of the ERL distribution parameters including: θ, λ, β, $a$ and $b$ which are scale and shape parameters. According to Cordeiro *et al*. (2011), the log likelihood function is given as:

$$L(\varphi) = n\log[B(a,b)] + \sum_{i=1}^{n}\log[K(y,\delta)] + (a-1)\sum_{i=1}^{n}\log[1+K(y,\delta)] + (b-1)[K(y,\delta)] \quad (27)$$

$$\varphi = (a,b,c,\delta) \text{ and } \delta = (\theta,\lambda,\beta) \text{ are vectors}$$

If c = 1, it becomes equation (27) which leads to

$$L(\varphi) = const - n\log[B(a,b)] + \sum_{i=1}^{n}\log[K(y,\delta)] + (a-1)\sum_{i=1}^{n}\log[1+K(y,\delta)] + (b-1)[K(y,\delta)] \quad (28)$$

$k(y,\delta)$ and $K(y,\delta)$ have been stated at the beginning.

The log likelihood function of ERL distribution is given as:

$$L_{BRLD}(\varphi) = -n\log[B(a,b)] + \sum_{i=1}^{n}[k(x)] + (a-1)\sum_{j=1}^{n}\log[K(x)] + (b-1)\sum_{i=1}^{n}\log[1-K(x)] \quad (29)$$

Taking the differentiation in respect to a, b, θ, λ and β give the following:

$$\frac{\partial L(\varphi)}{\partial a} = -n\frac{\Gamma'(a)}{\Gamma(a)} + n\frac{\Gamma'(a+b)}{\Gamma(a+b)} + \sum_{y=1}^{n}\log(1-K(Y,\delta)) \quad (30)$$

$$\frac{\partial L(\varphi)}{\partial b} = -n\frac{\Gamma'(b)}{\Gamma(b)} + n\frac{\Gamma'(a+b)}{\Gamma(a+b)} + \sum_{y=1}^{n}\log(K(Y,\delta)) \quad (31)$$

$$\frac{\partial L(\varphi)}{\partial \theta} = \sum_{j=1}^{n}\left[\frac{\frac{\partial[K(Y,\delta)]}{\partial \theta}}{K(Y,\delta)}\right] + (a-1)\sum_{y=1}^{n}\left[\frac{\frac{\partial[1-K(Y,\delta)]}{\partial \theta}}{1-K(Y,\delta)}\right] + (b-1)\sum_{y=1}^{n}\left[\frac{\frac{\partial[K(Y,\delta)]}{\partial \theta}}{(Y,\delta)}\right] \quad (32)$$

$$\frac{\partial L(\varphi)}{\partial \lambda} = \sum_{j=1}^{n}\left[\frac{\frac{\partial[K(Y,\delta)]}{\partial \lambda}}{K(Y,\delta)}\right] + (a-1)\sum_{y=1}^{n}\left[\frac{\frac{\partial[1-K(Y,\delta)]}{\partial \lambda}}{1-K(Y,\delta)}\right] + (b-1)\sum_{y=1}^{n}\left[\frac{\frac{\partial[K(Y,\delta)]}{\partial \lambda}}{(Y,\delta)}\right] \quad (33)$$

$$\frac{\partial L(\varphi)}{\partial \beta} = \sum_{j=1}^{n}\left[\frac{\frac{\partial[K(Y,\delta)]}{\partial \beta}}{K(Y,\delta)}\right] + (a-1)\sum_{y=1}^{n}\left[\frac{\frac{\partial[1-K(Y,\delta)]}{\partial \beta}}{1-K(Y,\delta)}\right] + (b-1)\sum_{y=1}^{n}\left[\frac{\frac{\partial[1-K(Y,\delta)]}{\partial \beta}}{(Y,\delta)}\right] \quad (34)$$

### 4.1 Analysis of Data

The data used for the analysis is a secondary data obtained from COVID-19 situation weekly epidemiological report 39; 5[th] – 11[th] July, 2021 (NCDC website state the website): Thirty-six (36) States including federal capital territory (FCT) with reported laboratory-confirmed COVID-19 cases, recoveries, deaths, samples tested and active cases (37 data points); and was accessed on Thursday 22[nd] July, 2021 put date accesses at reference not here. Only the death cases from all states of the federation are used for the analysis.

Table 1. Summary of Goodness of fit Statistics of Death Cases Data

| Minimum | Maximum | Skewness | Kurtosis | AD (P-value) | KS (p-value) | CVM (p-value) |
|---|---|---|---|---|---|---|
| 2.00 | 456.00 | 3.55 | 14.39 | 3.883e-12 | < 2.2e-16 | 3.591e-09 |

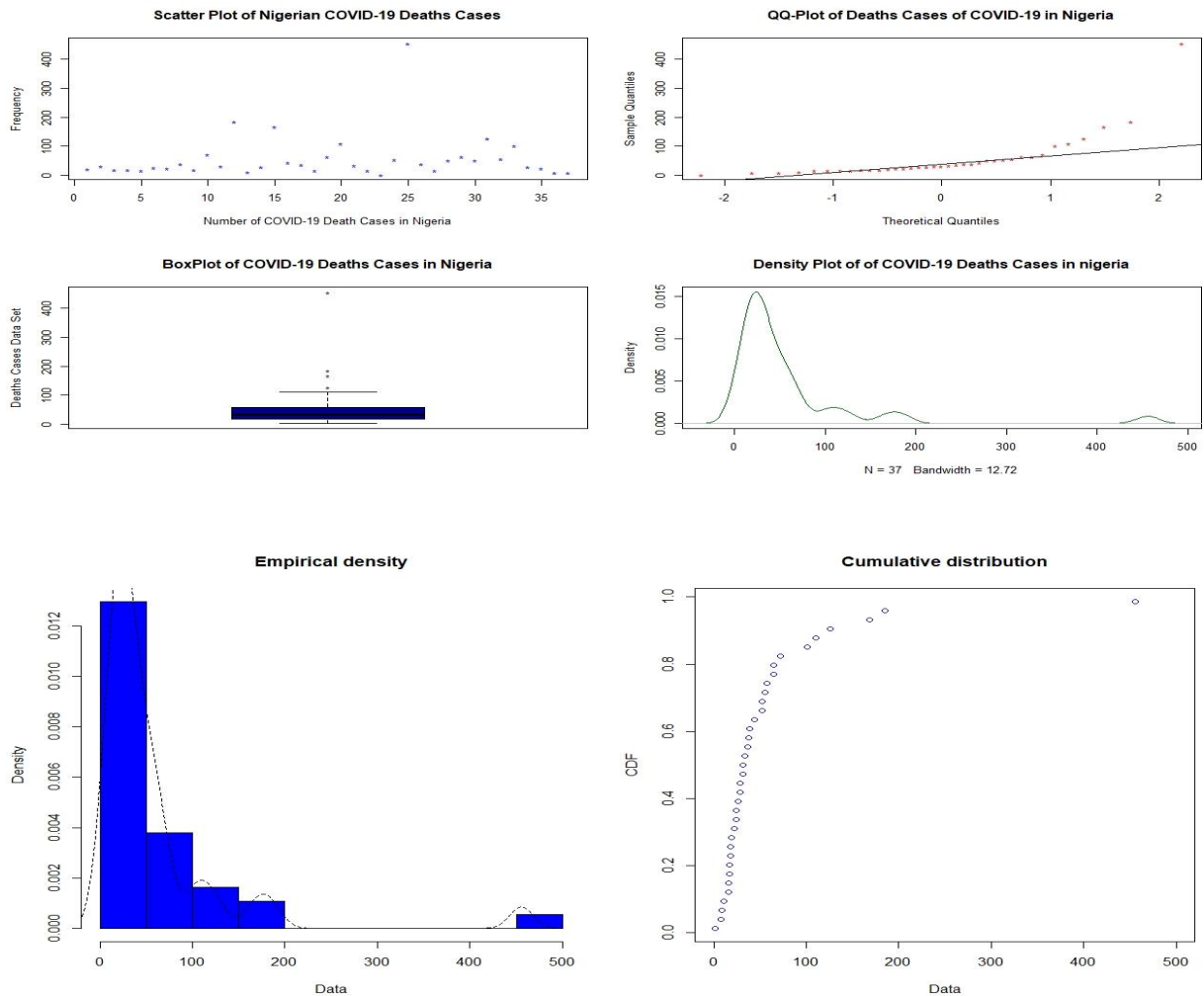

Figure 2: The Scatter, Theoretical Quantiles, Boxplot, Histogram, Density and distribution plot

Table 2: Contains the MLE, Standard Error (in parenthesis) and Model Selection Criteria.

| Model-Par /Model Sel | ERLD* | ExpLD | LRLD | BRD | RLD | ExpRLD | BLD |
|---|---|---|---|---|---|---|---|
| $a$ | **0.801 (0.015)** | 2.203 (0.044) | _ _ | 1.801 (0.038) | _ _ | 1.802 (0.040) | 0.801 (0.015) |
| $b$ | **5.500 (0.132)** | _ _ | 2.500 (0.056) | 0.999 (0.019) | _ _ | _ _ | 2.501 (0.057) |
| $\theta$ | **0.025 (NA)** | 0.002 (NA) | 0.008 (NA) | 0.048 (NA) | 0.000 (NA) | 0.016 (NA) | 0.010 (NA) |
| $\lambda$ | **0.113 (0.003)** | 0.109 (0.002) | 0.117 (0.003) | _ _ | 0.116 (0.003) | 0.110 (0.002) | 0.117 (0.003) |
| $\beta$ | **0.501 (0.014)** | _ _ | 0.500 (0.014) | 2.500 (0.063) | 2.900 (0.078) | 1.500 (0.037) | _ _ |
| $-2LogL$ | **1298.939** | 1629.850 | 2003.402 | 2201.424 | 2391.883 | 2498.683 | 2775.822 |
| $AIC$ | **2607.878** | 3265.700 | 4014.804 | 4410.848 | 4789.766 | 5005.366 | 5559.644 |
| $CAIC$ | **2609.813** | 3266.427 | 4016.054 | 4412.098 | 4790.493 | 5006.616 | 5560.916 |
| $HQIC$ | **2610.718** | 3267.404 | 4017.076 | 4413.120 | 4791.470 | 5007.638 | 5561.916 |
| $BIC$ | **2615.933** | 3270.533 | 4021.248 | 4417.292 | 4794.599 | 5011.810 | 5566.088 |

## 4.2  Result and Discussion

The summary of goodness of fit statistics is used to check for normality of the data; skewness, kurtosis, Anderson Darling (AD), Kolmogorov Smirnov (KS) and Cramer-Von-Mises (CVM) shown in Table 1 with their values clearly indicate that the data does not follow normal distribution since p-values less than 5%, skewness greater than 0 (zero) and kurtosis also greater than 3 (Karadimitriou and Shivam Mishra (2020). While, graphs from figure 2 show the nature of the data, the scatter, theoretical quantiles, boxplot, histogram, density and empirical cumulative distribution function (ecdf) plot show the data is skewed. For instance, non-linearity by scatter and quantiles plots, outliers by boxplot and skewness by histogram and density plots. The minimum and maximum values in the data set are inclusive.

The results obtained in Table 2 are based on parameter estimates by method of maximum likelihood estimation (MLEs). The standard error values are in bracket for all the models. The model ERLD is compared with other six models ExpLD, LRLD, BRD, RLD, ExpRLD and BLD Also, model selection criterion is performed on all models considered in the study. From the results, ERLD has the smallest values in all as we can see bold and starred where $\boldsymbol{AIC = 2607.878, CAIC = 2609.813, HQIC = 2610.718}$ and $\boldsymbol{BIC = 2615.933}$, which indicates that it is a robust and flexible model.

5.0   Conclusion

Despite the level of Nigerian COVID-19 death cases data set, the ERL distribution follows the movement of the data and has better representation of the data than any of the other existing distributions. The proposed distribution being flexible and versatile, can accommodate increasing, decreasing, bathtub and unimodal shape hazard function. It is therefore useful and effective in the analysis of clinical and survival data.

References

[1] N. I. Badmus, O. Faweya O and A. T. Adeniran. Modeling COVID-19 Pandemic Data with Beta Double Exponential Distribution. *Asian Journal of Research and Infectious diseases*. 5(4): (2020), 66 – 79. 63547. DOI: 10.9734/AJRID/2020/v5i430181. ISSN: 2582-3221.

[2] G. M. Cordeiro, and M. de Castro. A new family of generalized distributions. *Journal of Statistical Computation and Simulation*, 81(7), (2011), 883-898. https://doi.org/10.1080/00949650903530745

[3] A. H. El-Bassiouny, N. F. Abdoand and H. S. Shahen. Exponential Lomax Distribution. International Journal of Computer Applications. 121(13), (2015), 24-29.

[4] J. R. M. Hosking. L-moments Analysis and Estimation of distributions using Linear Combinations of Order Statistics. *Journal Royal Statistical Society B,* vol. 52, (1990), 105-124 http://www.jstor.org/stable/2345653

[5] M. C. Jones. Families of distributions arising from distributions of order statistics test. (2004), 13:1-43. DOI: https://doi.org/10.1007/BF02602999

[6] S. F. Karadimitriou. Checking normality in R – University of Sheffield. (2018) 1-4. /file/stcp https://vdocuments.mx/checking-normality-in-r

[7] Kawsar, Fatima., Uzma Jan. and S. P. Ahmad. Statistical Properties of Rayleigh Lomax distribution with applications in Survival Analysis. *Journal of Data Science*, (2018), 531-548. DOI: 10.6339/JDS.201807_16(3).0005

[8] NCDC (2021). General Fact Sheet-Data as at 11[th] July, 2021. States with Reported Laboratory Confirmed COVID-19 Cases, Recoveries, Deaths Samples Tested and Active Cases Accessed, July 22[nd], (2021)
.

[9] Shivam Mishra. Normality Test with Python in Data Science. (2020). Analytics Vidhya. https://medium.com/@shivamrkom?source=post_page........5abbefc81fdo

[10] M. M. Siddiqui. Some problems connected with Rayleigh distributions. J. Res. Nat. Bur. Stand, 60D, (1962), 167–174.

[11] K. S. Dey, M. M. Rahman, U. R. Siddiqi and A. Howlader. Analyzing the epidemiological outbreak of COVID19: A visual exploratory data analysis approach. Journal of Medical Virology.;92(6), (2020), 632-638. Available:https://doi.org/10.1002/jmv.25743

[12] World Health Organization. Surveillance case definitions for human infection with novel coronavirus (nCoV). https://www.who.int/publications-detail/surveillance-case-definitions-for-human-infection-with-novel-coronavirus-(ncov). Updated 2020. Accessed April 200, 2020.

[13] J. H. Yoo The fight against the nCoV outbreak: An arduous march has just begun. Journal of Korean Medical Science. (2020), 35(56). eISSN 1598 - 6357. ISSN 1011 - 8934. Accessed 22nd April, 2020.Available:https://doi.org/10.3346/jkms.2020.35.e56